\documentclass[runningheads]{llncs}

\usepackage{amsmath,amssymb,amsthm,color,enumerate}
\usepackage[normalem]{ulem}
\usepackage{mathtools}
\usepackage[sort&compress,comma,numbers]{natbib}
\usepackage[usenames,dvipsnames]{xcolor}
\usepackage[colorlinks=true, allcolors=blue]{hyperref}
\usepackage[nameinlink]{cleveref}
\usepackage[linesnumbered,ruled,vlined]{algorithm2e}
\usepackage{svg}
\usepackage{graphicx}
\usepackage{soul}
\usepackage{float}
\floatstyle{plaintop}
\restylefloat{table}

\graphicspath{./images/}

\theoremstyle{definition}
\newtheorem{assume}{Assumption}
\newtheorem{ex}{Example}

\newtheorem{conj}[ex]{Conjecture}

\crefname{assume}{Assumption}{Assumptions}
\crefname{equation}{}{}
\crefname{algorithm}{Algorithm}{Algorithms}
\crefname{figure}{Figure}{Figures}
\crefname{appendix}{Appendix}{Appendices}
\crefname{section}{Section}{Sections}
\crefname{table}{Table}{Tables}
\crefname{conj}{Conjecture}{Conjectures}

\newcommand{\R}{\mathbb{R}}

\newcommand{\Z}{\mathbb{Z}}

\newcommand{\optQ}{q^{\star}}
\newcommand{\optH}{h^{\star}}
\newcommand{\optT}{T^{\star}}
\newcommand{\RS}{\xi^{\star}}
\newcommand{\RSpred}{\hat{\xi^{\star}}}

\newcommand{\x}{\xi}
\allowdisplaybreaks
\begin{document}
\title{A learning-based algorithm to quickly compute good primal solutions for Stochastic Integer Programs}

\titlerunning{Learning based SIP}
%
\author{Yoshua Bengio\inst{2, 3} \and 
Emma Frejinger\inst{2,3} \and 
Andrea Lodi\inst{1, 3} \and 
Rahul Patel\inst{1, 3} \and 
Sriram Sankaranarayanan\inst{1}
%
}
\authorrunning{Y. Bengio et al.}
%
\institute{Canada Excellence Research Chair, Polytechnique Montreal \and
Department of Computer Science and Operations Research,
University of Montreal \and
Mila - Quebec Artificial Intelligence Institute\\
\email{andrea.lodi@polymtl.ca}}
\maketitle              
\begin{abstract}
We propose a novel approach using supervised learning to obtain near-optimal primal solutions for two-stage stochastic integer programming (2SIP) problems with constraints in the first and second stages. The goal of the algorithm is to predict a {\em representative scenario} (RS) for the problem such that, deterministically solving the 2SIP with the random realization equal to the RS, gives a near-optimal solution to the original 2SIP. Predicting an RS, instead of directly predicting a solution ensures first-stage feasibility of the solution. If the problem is known to have complete recourse, second-stage feasibility is also guaranteed. For computational testing, we learn to find an RS for a two-stage stochastic facility location problem with integer variables and linear constraints in both stages and consistently provide near-optimal solutions. Our computing times are very competitive with those of general-purpose integer programming solvers to achieve a similar solution quality.

\keywords{Stochastic Integer Programming  \and Machine Learning \and Heuristics.}
\end{abstract}
%
%
%

\section{Introduction}
Two-stage stochastic integer programming (2SIP) is a standard framework to model decision making under uncertainty. 
In this framework, first the so-called \emph{first-stage} decisions are made. 
Then, the values of some uncertain parameters in the problem are determined, as if sampled from a known distribution. 
Finally, the second set of decisions are made depending upon the realized value of the uncertain parameters, the so-called {\em second-stage} or {\em recourse} of the problem.
The decision maker, in this setting, minimizes the sum of (i) a linear function of the first-stage decision variables and (ii) the expected value of the second-stage optimization problem.

2SIP is studied extensively in the literature \citep{birge2011introduction,lulli2004branch,louveaux1992dual,kall1994stochastic, prekopa2013stochastic, shapiro2009lectures, dupavcova2003scenario, nemirovski2009robust, linderoth2006empirical,powell2015tutoriala, powell2015tutorialb} owing to its applicability in various decision making situations with uncertainty, like the stochastic unit-commitment problems for electricity generation \citep{powell2015tutoriala, powell2015tutorialb}, stochastic facility location problems \citep{louveaux1992dual}, stochastic supply chain network design \cite{santoso2005stochastic}, among others.
With the overwhelming importance of 2SIP a wide range of solution algorithms have been proposed, for example, \cite{ahmed2004finite, lulli2004branch, sen2005c, ahmed2013scenario, sen2010stochastic, caroe1998shaped}.

In this paper we are interested in using machine learning (ML) to obtain good primal solutions to 2SIP.
Along this line, \citet{nair2018} proposed a reinforcement learning-based {heuristic solver} to quickly find solutions to 2SIP. Given that the agent can be trained offline, the algorithm provided solutions much faster for some classes of problems 
compared to an open-source general-purpose mixed-integer programming (MIP) solver, in their case, SCIP \citep{GleixnerEtal2018OO, GleixnerEtal2018ZR}. However, their method is based on the following restrictive assumptions:
\begin{enumerate}[a.]
    \item \label{assum:allBin} All first-stage variables are required to be binary. General integer variables or continuous variables in the first stage cannot be handled. 
    \item \label{assum:allFeas}{\em Every} binary vector is required to be feasible for the first stage of the problem, i.e., no constraints are allowed in the first stage. 
\end{enumerate}
Assumption \ref{assum:allBin} above is intrinsic to the method in \citep{nair2018}, as both the {\em initialization policy} and {\em the local move policy} of the method involves flipping the bits of the first-stage decision vector. 
Hence, one cannot easily have general integer variables or continuous variables. 
Assumption \ref{assum:allFeas} is again crucial to the algorithm in \citep{nair2018}, as flipping a bit in the first stage could potentially make the new decision infeasible and it might require a more complicated feedback mechanism to check and discard infeasible solutions. 
In fact, if there are constraints, it is $NP$-hard to decide if there exists a flip that keeps the decision feasible. 
Alternatively, one could empirically penalize the infeasible solutions, but tuning the penalty might be a hard problem in itself. 

In contrast, our method does not require either of these two restrictive assumptions. We allow binary, general integer as well as continuous variables in both first and second stage of the problem. We also allow constraints in both stages of the problem.  Furthermore, we have a simple and direct approach to handle the first-stage constraints, without requiring any empirical penalties. 

We make the following common assumption to exclude pathological cases, where an uncertain realization can turn a feasible first-stage decision infeasible. 
\begin{assume}
The 2SIP has {\em complete recourse}, i.e., if a first-stage decision is feasible given the first-stage constraints, then it is feasible for all the second stage problems as well. \label{assum:compRec}
\end{assume}
We make another assumption of uncertainty with finite support, so we can have a proper benchmark to compare our solution against. However, this assumption can be readily removed, without affecting the proposed algorithm.
\begin{assume}
The uncertainty distribution in the 2SIP has a finite support. 
\label{assum:finiteSupport}
\end{assume}

\section{Problem definition}
We formally define a 2SIP as follows: 
\begin{subequations}
    \begin{align}
      \min_{x\in\R^{n_1}} & \quad c^Tx + \mathbb{E}_\xi \left[Q(x,\xi)\right]
      \\
      \text{subject to} & \quad Ax \leq b\\
      & \quad x_i \in \Z, \quad \forall \,i \in \mathcal{I}_1
    \end{align}
    \label{eq:2SIP}
\end{subequations}
where,
\begin{align*}
    Q(x,\xi) \quad &= \quad \min_{y\in\R^{n_2}} \left\{  q_\xi^Ty_{\xi}: Wy_{\xi} \leq h_{\xi} - T_{\xi}x, y_{\xi} \geq 0; y_i\in\Z\;\forall \,i\in \mathcal{I}_2 \right\}
\end{align*}
where $x  \in  \R^{n_1}$ and $y\in\R^{n_2}$ are the { first and second-stage decisions respectively}, $c  \in  \R^{n_1}$, $A  \in  \R^{m_1 \times n_1}$, $b  \in  \R^{m_1}$, $y_{\xi}  \in  \R^{n_2}$, $q_{\xi}  \in  \R^{n_2}$, $W  \in  \R^{m_2 \times n_2}$, $T_{\xi}  \in  \R^{m_2 \times n_1}$, $h_{\xi}  \in  \R^{m_2}$, $\mathcal{I}_1\subseteq \{1,\ldots,n_1\}, \mathcal{I}_2 \subseteq \{1,\ldots,n_2\}$.

When \cref{assum:finiteSupport} holds, the 2SIP described above can also be expressed as a single deterministic MIP as follows:
\begin{subequations}
    \begin{align}
      \min_{x, y} & \quad c^Tx + \sum_{\forall \xi \in \Xi} p_{\xi} q_{\xi}^T y_{\xi} 
      \\
      \text{subject to} & \quad Ax \leq b
      \\
      & \quad Wy_{\xi} \leq h_{\xi} - T_{\xi} x \qquad \forall \xi \in \Xi 
      \\
      & \quad x_i \in \Z, \quad \forall \,i \in \mathcal{I}_1
      \\
      & \quad y_{\xi i} \in \Z, \quad \forall \xi \in \Xi, \forall \,i \in \mathcal{I}_2.
    \end{align}\label{eq:DetermSIP}
\end{subequations}

When \cref{assum:finiteSupport} does not hold, the formulation \cref{eq:DetermSIP} could be a finite-sample approximation of \cref{eq:2SIP}, which is extensively studied in the stochastic programming literature. Imitating \cite{nair2018}, we compare our algorithm against solving \cref{eq:DetermSIP} with a general-purpose MIP solver.



\section{Methodology}

In this section, we discuss the algorithmic contribution of the paper.
\subsection{Surrogate formulation}\label{sec:Surr}
We first define the objective value function (OVF)   $\Phi:\R^{n_1}\to\R$, mapping $x\mapsto c^Tx + \mathbb{E}_\xi \left[Q(x,\xi)\right]$ - the function we are trying to optimize over the mixed-integer set defined in \cref{eq:2SIP}. 

Given \cref{eq:DetermSIP}, we define {\em the surrogate problem associated with $\bar \xi = (\bar q,\bar h, \bar T)$}, as follows:
\begin{subequations}
    \begin{align}
        \min_{x,y} & \quad c^Tx + \bar q^{ T} y
      \\
      \text{subject to} & \quad Ax \leq b
      \\
      & \quad Wy \leq \bar h- \bar Tx
      \\
      & \quad x_i, y_j \in \Z, \quad \forall \,i \in \mathcal{I}_1;\,j \in \mathcal{I}_2
    \end{align} \label{eq:RLrefor}
\end{subequations}

In other words, should the value that the uncertain parameters are going to take is deterministically known to be $\bar \xi$, then the decision maker can solve the surrogate problem associated with $\bar\xi$.
Now, the idea behind the algorithm proposed in this paper is captured by \cref{conj:xiHatExists}. 
\begin{conj} 
\label{conj:xiHatExists}
Let \cref{eq:DetermSIP} (and hence \cref{eq:2SIP}) have an optimal objective value of $f^*$. There exists $\optQ, \optH, \optT$ in $\R^{n_2}$, $\R^{m_2}$ and $\R^{m_2\times n_1}$ such that if $(x^\dagger, y^\dagger)$ solves the (much smaller) surrogate problem defined by $(\optQ, \optH, \optT)$, then, $ f^* = \Phi(x^\dagger)$.
\end{conj}

Observe that by construction, $x^\dagger$ is feasible to the original problem in \cref{eq:2SIP}. Also, \cref{conj:xiHatExists} asserts that, there exists a realization of the uncertainty ($\RS = (\optQ, \optH, \optT)$) such that if one deterministically optimizes for that realization $\RS$, then its solutions are optimal for the original 2SIP. Each such $\RS$ is called a {\em representative scenario} (RS) for the given 2SIP. 

Now, given adequate computing resources, one can solve the following bilevel program to obtain an RS. \begin{subequations}
\begin{align}
    \min_{\substack{U,v,w\\x,y}} \quad&\quad c^Tx +  \sum_{\forall \xi \in \Xi} p_{\xi} q_{\xi}^T y_{\xi} 
      \\
      \text{subject to} & 
      \quad  (x,w) \in \arg\min_{x,w} 
      \left \{c^Tx + v^Tw:
      \begin{array}{cc}
      Ax \leq b; & \\
      Ww \leq v - Ux; &\\
      x_i \in \Z     &  \forall \,i \in \mathcal{I}_1\\
      w_i \in \Z       & \forall \,i \in \mathcal{I}_2
      \end{array}
      \right \}
      \\
      & \quad Wy_{\xi} \leq h_{\xi} - T_{\xi} x \qquad \forall \xi \in \Xi 
      \\
      & \quad y_{\xi i} \in \Z, \quad \forall \xi \in \Xi, \forall \,i \in \mathcal{I}_2
\end{align}\label{eq:bilevel}
\end{subequations}
If the optimal value of this problem matches the optimal value of the original 2SIP, then the corresponding values for $(U, v, w)$ form the RS.
Note that if $T_\xi = T,\,\forall \xi \in \Xi$, then \cref{eq:bilevel} is a mixed-integer bilevel {\em linear} program (MIBLP) and can hopefully be solved faster than the general case. 





\subsection{Learning algorithm}

    The goal of ML algorithms is to predict an optimal ($U,v,w$) to \cref{eq:bilevel}, given the data for the 2SIP. On the one hand, we are expecting the ML algorithms to predict the solutions of a seemingly much harder optimization problem than the original 2SIP. On the other hand, this is easier for ML since there are no constraints on the predicted variables -- $U,v,w$.
    Supervised learning is a natural tool to achieve this goal.

{Supervised learning} can be used if there is a training dataset of problem instances and their corresponding RS. The task of predicting RS can be formulated as a regression task as RS is real valued. The algorithm tries to minimize the mean squared error (MSE) between the true and predicted RS. The prediction can also be evaluated on the merits of optimization metrics, comparing the solution and objective value of true and predicted RS. 

%

%
\section{Computational study}\label{sec:study} 
This section discusses the computation study performed to support \cref{conj:xiHatExists}. 
%
\subsection{Problem definition} 

In this work, we consider a version of two-stage stochastic capacitated facility location (S-CFLP) for computational analysis. The problem is enhanced such that both the first and the second stage of the problem have integer as well as continuous variables. 
More precisely, the first stage consists of deciding (i) the locations where a facility has to be opened (binary decisions), and (ii) if a facility is open, then the maximum demand that the facility can serve (continuous decisions). There are also constraints which dictate bounds on the total number of facilities that can be opened. 
The uncertainty in the problem pertains to the demand values at various locations in the second stage of the problem, which are sampled from a finite distribution. Once the demand is realized, the second stage consists in deciding (i) if a given open facility should serve the demand in a location (binary decisions) (ii) if a facility serves the demand in a location, then what quantity of demand is to be served (continuous decisions). These decisions have to ensure that the demand and supply constraints are met. The problem formulation is presented formally in \cref{app:probfor}.



\subsection{Data generation}
\paragraph{Generate instances.}
We generate 50K instances of S-CFLP, with 10 facilities, 10 clients and 50 scenarios. We provide details on how the data for these instances are generated in \cref{app:DatGen}. The generated instances are solved using Gurobi, running 2 threads, to at most 2\% gap or 10~min time limit. 

Next, we compute an RS for each of the 50K instances. As stated earlier, one could solve mixed-integer bilevel programs \cref{eq:bilevel} to obtain the RSs. However, due to the computational burden, we use heuristics that work using the knowledge of the (nearly) optimal solutions to the 2SIP already obtained from Gurobi. These heuristics are detailed in \cref{app:heuristics}. Out of 50K instances, they find an RS for 49,290 instances. We believe that a more thorough search will enable us to find the RS for all the problems. 

\subsection{Learning algorithm}
We formulate the task of predicting the RS as a regression task. The size of the dataset, which comprises of instances and their corresponding RS, is 49,290. The dataset is split into training and test sets of size 45K and 4,290, respectively. Further, a validation set of 5K is carved out from the training set.
We use linear regression (LR) and artificial neural network (ANN) to minimize the MSE between the true and predicted RS. 

\paragraph{Feature engineering.}
It is well known that features describing the connection between variables, constraints and other interaction help ML to perform well rather than just providing plain data matrices \citep{khalil2016learning, bonami2018learning, gasse2019exact, bengio2018machine}. In this spirit, along with the fixed and variable costs to open facilities at different locations, we also provide aggregated features on the set of scenarios. These features give information about each of the potential locations for facilities in S-CFLP as well as the way different locations interact through the demands in adjacent nodes. A detailed set of the features is given in \cref{app:feature_engineering}.

\subsection{Comparison}
In order to evaluate the ML-based prediction of $\RS$, which we refer to as $\RSpred$, we compare 
the solution obtained by solving the surrogate problem associated with $\RSpred$ against solutions obtained by various algorithms.

We use LR and ANN to predict $\RS$.
We compare these predictions against (i) Solution obtained using Gurobi by solving \cref{eq:DetermSIP} (GRB) (ii) a solution obtained by solving the surrogate associated with the average scenario, namely $\sum_{i=1}^N\Xi_i/N$ (AVG) (iii) a solution obtained by solving the surrogate associated with a randomly chosen scenario from the $N$ choices (RND) (iv) a solution obtained by solving the surrogate associated with a randomly chosen scenario from the distribution of the scenario predicted by LR (DIST). Note that GRB produces better solutions (in most cases) than the ML methods, however, taking a significantly longer time. We therefore assess the time it takes GRB to get a solution of comparable quality to LR and ANN. We refer  to these as GRB-L and GRB-A, respectively.

\section{Results}
\begin{figure}[!htbp]
    \centering
    \includegraphics[width=12cm]{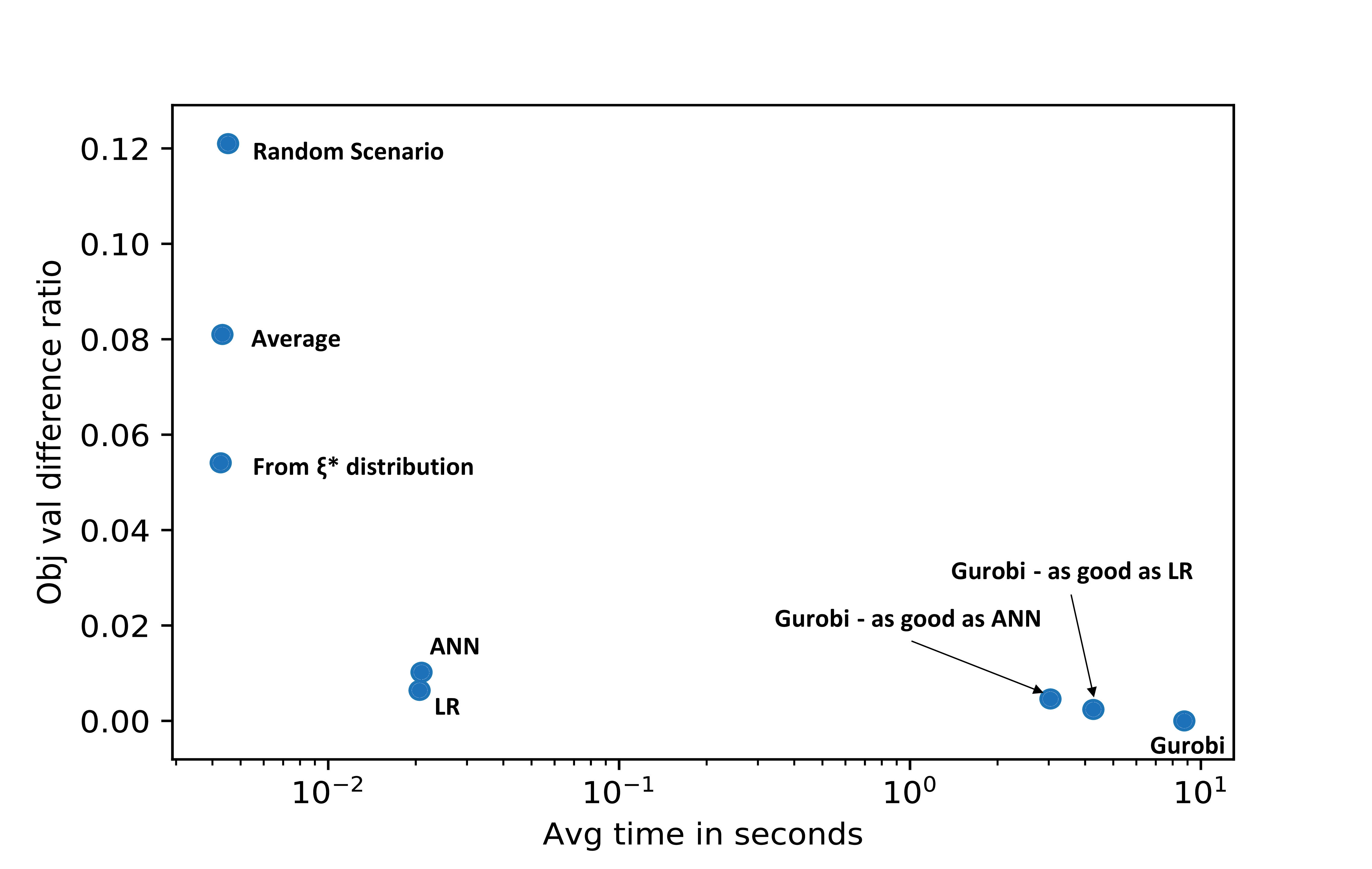}
    \caption{Objective value difference ratio vs. avg time in seconds}
    \label{fig:obj_vs_time}
\end{figure}

\cref{tab:obj_val_diff} reports the {\em objective value difference ratio} defined as $(($Obj val by a method$- $GRB obj val$)/$GRB obj val$)$ for each method and \cref{tab:time_stats} statistics on computing times. Before analyzing the results in more detail, we note a key finding that emerges. Namely, LR and ANN perform almost as good as GRB (in terms of quality of the objective value) in a fraction of the time taken by GRB. \cref{fig:obj_vs_time} captures the trade off between the quality of solutions obtained by different methods as well as the time taken to obtain these solutions.

We observe from \cref{tab:obj_val_diff} that LR and ANN produce decisions that are as good as GRB ones on an average (and by the median value), and in some cases the ML-based methods even beat GRB, i.e., produce solutions whose objective value is better than that of GRB. This is possible because GRB does not necessarily solve the problem to optimality, but only up to a gap of 2\%. Further, even in the worst of the 4,290 test cases, LR is at most 2.64\% away from GRB. To show that this is not easily achieved, we also compare GRB against AVG, RND and DIST. We observe that these methods perform much poorer than GRB, unlike LR and ANN.


\begin{table}[!htbt]
    \setlength{\tabcolsep}{5pt}
    \centering
    \begin{tabular}{|r||r r r || r r|}
    \hline
    \textbf{GRB vs.} & \textbf{AVG} & \textbf{RND} & \textbf{DIST} & \textbf{LR} & \textbf{ANN}\\
    \hline
    \hline
    Min  & 3.36  & -0.33  & -0.17    &  -0.62   & -0.45 \\
    Max  & 14.23  &  94.42  & 49.87  &   2.64 & 7.85 \\
    Avg.  & 8.10& 12.10  & 5.41  &  0.64 & 1.02 \\
    Median   & 8.08 & 8.24  & 3.54 &   0.60  & 0.90 \\
    Std. dev. &   1.59 &  12.11  & 5.39 &  0.41 & 0.70 \\
    \hline
    \end{tabular}
    \caption{Objective value difference ratio, GRB vs. the five other methods (in \%)}
    \label{tab:obj_val_diff}
\end{table}

\begin{table}[!htbp]
    \setlength{\tabcolsep}{5pt}
    \centering
    \begin{tabular}{|r||r c c c c c c c|}
    \hline
    \textbf{Method} & \textbf{GRB} & \textbf{AVG} & \textbf{RND} & \textbf{DIST} & \textbf{GRB-L} & \textbf{GRB-A} & \textbf{LR} & \textbf{ANN} \\ \hline
    Min  & 0.4354 & 0.0041 & 0.0041 & 0.0040 & 0.4049  & 0.4268 & 0.0194 & 0.0197   \\ \hline
    Max  & 559.32 & 0.0077 & 0.0081& 0.0065 & 140.43 & 140.43 & 0.0454 & 0.0457   \\ \hline
    Avg.  & 8.7713 & 0.0043 & 0.0045& 0.0043 & 4.2650 & 3.0398 & 0.0206 & 0.0209   \\ \hline
    Median  & 2.2621 & 0.0043 & 0.0044& 0.0042 & 1.4613 & 1.2898 & 0.0200 & 0.0204   \\ \hline
    Std. dev. & 17.919 & 0.0001 & 0.0003& 0.0001 & 8.5576 & 6.8362 & 0.0019 & 0.0019   \\ \hline
    \end{tabular}
    \caption{Statistics on computing times of the different methods}
    \label{tab:time_stats}
\end{table}

Analyzing the time improvement in using LR and ANN, we observe from \cref{tab:time_stats} that these methods solve the S-CFLP {\em orders of magnitude} faster than GRB. 
Indeed, GRB takes over 8 seconds on an average to solve these problems, while the maximum time is 0.046 seconds using LR and ANN. We emphasize that the time taken to solve using the ML methods {\em includes} the time elapsed in computing the values of the features used in ML and the time elapsed in solving the surrogate associated with the corresponding $\RSpred$.
Recall that GRB-L and GRB-A denote the time it takes GRB to produce a solution of comparable quality to LR and ANN. The results show that GRB cannot produce a solution of the same quality as LR and ANN in a comparable time. In fact, LR and ANN are still orders of magnitude faster than GRB.



\section{Discussion} 

In this paper, we present an algorithm to solve 2SIP using ML-based methods. The method hinges on the existence of the RS conjectured in \cref{sec:Surr}. Computationally, we see that the methods proposed in this paper consistently provide good quality solutions to  S-CFLP orders of magnitude faster.

An important observation we had while training the models is that we were never able to get the training loss close to zero. Naturally, the predicted RS in the test dataset were quite different from the RS estimated using our heuristics. The differences in the predicted values of the components of RS and those obtained using the heuristics are documented by \cref{fig:hist} in the Appendix. Despite this, the solution value to the 2SIP as determined by our algorithm were near optimal as shown in the results and significantly better than those obtained with other methods. The mismatch between the ML metrics and those characterizing discrete optimization problems is a known issue \citep{bengio2018machine} requiring extensive research and, in our context, we believe that exploring this avenue might produce better solutions.

Another interesting observation is that LR beats ANN in this task. We suspect that this is partly caused by the parsimony offered by LR. However, this is also encouraging news that the sample complexity of the learning task might be relatively small in general. We believe that a natural extension to this work is to provide these analyses more formally. 

Further, we believe that computational tests assessing the performance of the algorithms in different datasets of 2SIP is crucial to show how much and where our method generalizes. This might involve learning solutions to other forms of 2SIP like the stochastic unit-commitment problem, the stochastic supply chain-network design problem etc. These are cases where we believe that \cref{conj:xiHatExists} still holds, but we do not have computational validation. 

Finally, we would also be interested in extending the theory side when \cref{conj:xiHatExists} is not expected to hold at all or holds only with weaker guarantees; for example, in the case where both the stages are mixed-integer nonlinear programming (MINLP) problems. In such cases, it will be useful to understand the reach of ML-based solution techniques as opposed to traditional MINLP solvers.

\bibliographystyle{plainnat}
\bibliography{references}

\appendix
\section {Computational test details}
\subsection{Problem Formulation}
\label{app:probfor}
We provide below the problem considered in this work for computational study.
\begin{subequations}
\begin{align}
\min_{b\in\{0,1\}^n, v\in\R^n_+} 
\quad & \quad 
\sum_{i=1}^n 
\left (c^f_i b_i + c^v_iv_i\right ) + \mathbb{E}_{\xi}(Q(x,\xi))\label{eq:CFL:obj}\\
& \quad \frac{n}{10} \quad\leq\quad \sum_{i=1}^n b_i \quad\leq\quad \frac{3n}{4} \label{eq:CFL:bnd}\\
& 
v_i \quad\leq\quad Mb_i \label{eq:CFL:bM}\\
\text{where } Q(x,\xi)\text{ is } &\nonumber\\
\min_{u\in\{0,1\}^{n\times n}, y\in\R^{n\times n}_+} \quad & \quad \sum_{i=1}^n\sum_{j=1}^n c^{tv}_{ij}y_{ij} + \sum_{i=1}^n\sum_{j=1}^n c^{tf}_{ij}u_{ij} \label{eq:CFL:2obj}\\
& \quad \sum_{j=1}^n y_{ij} \quad \leq\quad v_i \label{eq:CFL:suppl}\\
& \quad \sum_{i=1}^n y_{ij} \quad =\quad d_j(\xi)\label{eq:CFL:dem}\\
& \quad y_{ij} \quad \leq\quad u_{ij}M \label{eq:CFL:uM}
\end{align}\label{eq:CFL}
\end{subequations}
In this problem, we minimize the fixed and variable costs of opening a facility along with the fixed and variable costs of transportation between the facilities and the demand nodes. There are $n$ potential locations where a facility could be opened. A fixed cost of $c^f_i$ is incurred, if a facility is opened in location $i$, and a variable cost of $c^v_i$ is incurred per-unit capacity of the facility opened in location $i$. The binary variable, $b_i$ tracks if a facility is opened in location $i$ and the continuous variable $v_i$ indicates the size of the facility at location $i$. The constraint in \cref{eq:CFL:bM}, along with the binary constraint on $b$ ensures that the costs are incurred in the right way. Finally, \cref{eq:CFL:bnd} is a complicating constraint, which says that at least a tenth of the locations must have a facility open and not more than three-quarters of the locations must have a facility open. 

In the second stage, $c^{tf}_{ij}$ is the fixed cost incurred in transporting from location $i$ to $j$; $c^{tv}_{ij}$ is the per-unit variable cost incurred in transporting from $i$ to $j$. The binary variable $u_{ij}$ denotes if any transportation happens from $i$ to $j$ and the continuous variable $y_{ij}$ denotes the actual quantity transported from $i$ to $j$.
The second-stage objective in \cref{eq:CFL:2obj} minimizes the transportation cost incurred under a random demand scenario parameterized by $\xi$.
Then, \cref{eq:CFL:suppl} ensures that the total quantity transported out of a facility is not greater than the capacity of the facility, while \cref{eq:CFL:dem} ensure that the total quantity supplied to a location $j$ equals the (random) demand at $j$. Finally, constraints \cref{eq:CFL:uM} link $u$ and $y$ variables appropriately.

\subsection{Data generation}
\label{app:DatGen}

\paragraph{Instance generation. }We generate 50K instances of S-CFLP, with 10 facilities, 10 clients and 50 scenarios. The random parameters $c^f$, $c^v$ and $\Xi=[\xi_1,\dots,\xi_{50}]$ vary across instances, where as $c^{tf}$ and $c^{tv}$ are fixed across all instances. Moreover, $c^f$ and $c^v$ are sampled from a discrete uniform distribution [15, 20) and [5, 10), respectively. The demand matrix $\Xi$ is generated by first evaluating $\lambda = \left\lfloor (c^f + 10*c^v)/\sqrt{n}\right\rfloor$. The $i^{th}$ demand scenario ($1 \leq i \leq 50$) is generated by sampling from a Poisson distribution with the mean equal to $\lambda$.

The generated instances are solved using Gurobi, running 2 threads, to optimality (less than 2\% gap) or 10~min time limit. We store the objective value, gap closed and master solution $x^* = (b^*$, $v^*)$. We were able to solve all the instances up to the specified gap, within the specified time limit.

We follow Algorithm~\ref{alg:generate_xi_hat} for generating $\RS$. Then, $|\Xi|$ refers to the cardinality of $\Xi$ and $c=1.01$ in step 6. 
\begin{algorithm}[t]
\KwData{2SIP $\mathcal{P}$ with objective value $o^e$ and first-stage solution $x^e$, max iterations $iter$}
\KwResult{$\RS$ or \textit{NULL}}
    $\bar\xi \leftarrow \frac{1}{|\Xi|} \sum_{i=1}^{|\Xi|} \xi_i$\;
    \While{iter}{
        Formulate surrogate problem $\mathcal{P}^{'}$ using $\bar\xi$\;
        Solve $\mathcal{P}^{'}$ and extract first-stage solution $x^{\bar\xi}$\; $o^{\bar\xi} = \Phi(\bar\xi)$
        \If{$o^{\bar\xi} \leq c * o^{e}$}{ 
            \Return $\bar\xi$\;
        }
        \Else{
            Perturb $\bar\xi$ using heuristics based on $x^e$ and $x^{\bar\xi}$
        }
    }
    \Return{NULL}
\caption{\textsc{GenerateXiHat}}
\label{alg:generate_xi_hat}
\end{algorithm}
The heuristics for updating the RS, based on $x^e$ and $x^{\bar\x}$, are described in \cref{app:heuristics}.
\subsection{Heuristics}
\label{app:heuristics}
Let $x^* = (b^*, v^*)$ and $x^{\bar\x} = (b^{\bar\x}, v^{\bar\x})$ be the first-stage optimal and surrogate solution associated with the scenario $\bar\x$, respectively. 
There are three heuristics that we use in tandem to generate $\RS$. 
The first heuristic is based on the comparison of facilities being open or close in the optimal and surrogate solution. 
The demand in the $\bar\x$ is zeroed out at clients for which the $b^*$ is close and $b^{\bar\x}$ is open, as suggested by 
\begin{equation}
    \label{eq:heuristic_zero_demand}
    b^*_i = 0 \land b^{\bar\x}_i = 1 \implies \bar\x_i = 0 \quad i=1,\dots,n. 
\end{equation}
The remaining two heuristics are based on the comparison of capacities installed in facilities in the optimal and the surrogate solution. 
First, the client with maximum absolute difference between capacities installed in optimal and surrogate solution is identified i.e., $\operatorname*{argmax}_i |v^*_i - v^{\bar\x}_i|$. 
The demand at this client in the $\bar\x$ is updated either by a fixed percentage $p$ of the current demand 
\begin{equation}
    \label{eq:heuristic_percent}
   \bar \x_i = \bar\x_i + \frac{v^*_i - v^{\bar\x}_i}{|v^*_i - v^{\bar\x}_i|} \times  p\bar\x_i,
\end{equation}
or by a fraction $f$ of the difference between the capacities installed in optimal and surrogate solutions
\begin{equation}
    \label{eq:heuristic_difference}
    \bar\x_i = \bar\x_i + (v^*_i - v^{\bar\x}_i)  f  \bar\x_i.
\end{equation}

\subsection{Feature engineering}
\label{app:feature_engineering}
\begin{figure}[t]
    \centering
    \includegraphics[height=6cm]{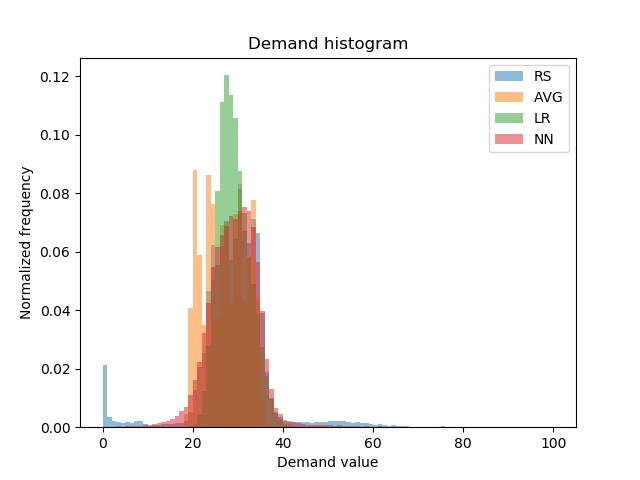}
    \caption{Demand histogram for different methods.}
    \label{fig:hist}
\end{figure}
The inputs of the models are $c^f, c^v$ and $\Xi$. We do not provide $c^{tf}$ and $c^{tv}$ in the input as they are fixed across all instances. We do feature engineering on $\Xi$, instead of providing it as a raw input, to extract information which can be useful in predicting $\RS$. Let $\Xi$ be an $m \times n$ matrix, where $m$ is the number of scenarios and $n$ is the number of clients. We calculate minimum, maximum, average, standard deviation, median, $75^{th}$ quantile, and $25^{th}$ quantile of $\Xi_{[:,i]}$ ($i^{th}$ column of $\Xi$) for $i=1, \dots, n$.

We also find the percentage of scenarios in which some fraction of demand for a client is greater than and less than the demand on all the other nodes 
$$
\frac{c * \Xi_{[:, i]} \geq \Xi_{[:, \neq i]}}{m}
\qquad\text{and}\qquad
\frac{c * \Xi_{[:, i]} \leq \Xi_{[:, \neq i]}}{m}.
$$
We set $c$ to different values (0.9, 1, 1.1, 1.2 and 1.5) and we thus end up with an input vector of size 190, combining $c^f$, $c^v$ and features extracted from $\Xi$.
The feature extraction performs an aggregation over the number of scenarios.




%




\end{document}